\documentclass{fpsac}

    \usepackage{amssymb, amsthm, amsfonts}
    \usepackage{graphicx}
    \usepackage{pstricks, pst-plot}

    \newtheorem{thm}{Theorem}[section]
    \newtheorem{prop}[thm]{Proposition}
    \newtheorem{lem}[thm]{Lemma}
    \newtheorem{cor}[thm]{Corollary}

    \theoremstyle{definition}
        \newtheorem{defn}[thm]{Definition}
        \newtheorem{exgr}[thm]{Example}

        \newtheorem{alg}[thm]{Algorithm}

    \theoremstyle{remark}
        \newtheorem{rem}[thm]{Remark}

    \numberwithin{equation}{section}

    \author[A. Burstein]{Alexander Burstein}
    \address{Department of Mathematics, Iowa State University, Ames, IA 50011-2064, USA}
    \email{burstein@math.iastate.edu}
    \urladdr{http://www.math.iastate.edu/burstein/}
    \thanks{The work of the first author was supported in part by the
    U.S. National Security Agency Young Investigator Grant H98230-06-1-0037.}

    \author[I. Lankham]{Isaiah Lankham}
    \address{Department of Mathematics, University of California, Davis, CA 95616-8633, USA}
    \email{issy@math.ucdavis.edu}
    \urladdr{http://www.math.ucdavis.edy/$\sim$issy/}
    \thanks{The work of the second author was supported in part by the
    U.S. National Science Foundation under Grants DMS-0135345 and
    DMS-0304414.}

    \title[Restricted Patience Sorting]{Restricted Patience Sorting
    and Barred Pattern Avoidance}



\begin{document}

    \begin{abstract}

       Patience Sorting is a combinatorial algorithm that can be
       viewed as an iterated, non-recursive form of the Schensted
       Insertion Algorithm.  In recent work the authors have shown
       that Patience Sorting provides an algorithmic description for
       permutations avoiding the barred (generalized) permutation
       pattern $3\mathrm{-}\bar{1}\mathrm{-}42$.  Motivated by this
       and a recently formulated geometric form for Patience Sorting
       in terms of certain intersecting lattice paths, we study the
       related themes of restricted input and avoidance of similar
       barred permutation patterns.  One such result is to
       characterize those permutations for which Patience Sorting is
       an invertible algorithm as the set of permutations
       simultaneously avoiding the barred patterns
       $3\mathrm{-}\bar{1}\mathrm{-}42$ and
       $3\mathrm{-}\bar{1}\mathrm{-}24$.  We then enumerate this
       avoidance set, which involves convolved Fibonacci numbers.

       \begin{resume}
       \emph{Patience Sorting} est un algorithme combinatoire que l'on
       peut comprendre comme \'{e}tant une version it\'{e}r\'{e}e,
       non-r\'{e}cursive de la correspondence de Schensted.  Dans leur
       travail r\'{e}cent les auteurs ont d\'{e}montr\'{e} que
       \emph{Patience Sorting} donne une description algorithmique des
       permutations \'{e}vitant le motif barr\'{e}
       (g\'{e}n\'{e}ralis\'{e}) $3\mathrm{-}\bar{1}\mathrm{-}42$.
       Motiv\'{e}s par ceci et par une forme r\'{e}cemment
       formul\'{e}e de \emph{Patience Sorting} en termes de certaine
       parcours du treillis intersectants, nous \'{e}tudions les
       th\`{e}mes connexe d'input restreinte et permutations qui
       \'{e}vitent de similaire motifs barr\'{e}s.  Un de nos
       r\'{e}sultats est de caract\'{e}riser les permutations pour
       lesquelles \emph{Patience Sorting} est un algorithme inversible
       comme \'{e}tant l'ensemble des permutations \'{e}vitant
       simultan\'{e}ment les motifs barr\'{e}s
       $3\mathrm{-}\bar{1}\mathrm{-}42$ and
       $3\mathrm{-}\bar{1}\mathrm{-}24$.  Nous \'{e}num\'{e}rons
       ensuite cet ensemble, qui utilise des convolutions des nombres
       de Fibonacci.
       \end{resume}

    \end{abstract}

    \maketitle  %%% make the title

    \section{Introduction}
    \label{sec:intro}

    The term \emph{Patience Sorting} was introduced in 1962 by C.~L.
    Mallows \cite{refMallows1962, refMallows1963} while studying a
    card sorting algorithm invented by A.~S.~C. Ross.  Given a
    shuffled deck of cards $\sigma = c_{1} c_{2} \cdots c_{n}$ (which
    we take to be a permutation $\sigma \in \mathfrak{S}_{n}$), Ross
    proposed the following algorithm:\medskip

    \begin{quote}
    \label{alg:PatienceSortingProcedure}

        \begin{enumerate}
            \item[\qquad{\bf Step 1}] Use what Mallows called a
            ``patience sorting procedure'' to form the subsequences
            $r_{1}, r_{2}, \ldots, r_{m}$ of $\sigma$ (called
            \emph{piles}) as follows:\medskip

                \begin{itemize}
                    \item Place the first card $c_{1}$ from the deck
                    into a pile $r_{1}$ by itself.\smallskip

                    \item For each remaining card $c_{i}$ ($i = 2,
                    \ldots, n$), consider the cards $d_{1}, d_{2},
                    \ldots, d_{k}$ atop the piles $r_{1}, r_{2}, \ldots,
                    r_{k}$ that have already been formed.\smallskip

                        \begin{itemize}
                            \item If $c_{i} > \max\{d_{1}, d_{2},
                            \ldots, d_{k}\}$, then put $c_{i}$ into a
                            new right-most pile $r_{k+1}$ by
                            itself.\medskip

                            \item Otherwise, find the left-most card
                            $d_{j}$ that is larger than $c_{i}$ and put
                            the card $c_{i}$ atop pile
                            $r_{j}$.\bigskip
                        \end{itemize}
                \end{itemize}

            \item[{\bf Step 2}] Gather the cards up one at a time from
            these piles in ascending order.\bigskip

        \end{enumerate}

    \end{quote}

    \setlength{\textheight}{8.5in}

    \noindent We call {\bf Step 1} of the above algorithm
    \emph{Patience Sorting} and denote by $R(\sigma) = \{r_{1}, r_{2},
    \ldots, r_{m}\}$ the \emph{pile configuration} associated to the
    permutation $\sigma \in \mathfrak{S}_{n}$.  Moreover, given any
    pile configuration $R$, one forms its \emph{reverse patience word}
    $RPW(R)$ by listing the piles in $R$ ``from bottom to top, left to
    right'' (i.e., by reversing the so-called ``far-eastern
    reading'').  In \cite{refBLFPSAC05} these words are characterized
    as being exactly the elements of the avoidance set
    $S_{n}(3\textrm{-}\bar{1}\textrm{-}42)$.  That is, they are
    permutation avoiding the generalized pattern $2\textrm{-}31$
    unless every occurrence of $2\textrm{-}31$ is contained within an
    occurrence of the generalized pattern $3\textrm{-}1\textrm{-}42$.
    (A review of generalized permutation patterns can be found in
    Section \ref{sec:Introduction:GeneralizedPatternAvoidance} below).

    We illustrate the formation of $R(\sigma)$ and $RPW(R)$ in the
    following example.

    \begin{exgr}\label{eg:NormalPSexample}
        Let $\sigma = 6 4 5 1 8 7 2 3 \in \mathfrak{S}_{8}$.  Then we
        form the pile configuration $R(\sigma)$ as follows:\\

        \begin{center}

            \begin{tabular}{l p{56pt} l p{56pt} l p{56pt}}
                \begin{minipage}[c]{68pt}
                    \begin{flushleft}
                    Form a new pile with \textbf{6}:
                    \end{flushleft}
                \end{minipage}
                &
                \begin{tabular}{l l l}
                                      & & \\
                                      & & \\
                       \textbf{6} & &
                \end{tabular}
                &
                \begin{minipage}[c]{68pt}
                    \begin{flushleft}
                        Then place \textbf{4} atop 6:
                    \end{flushleft}
                \end{minipage}
                &
                \begin{tabular}{l l l}
                                      & & \\
                       \textbf{4} & & \\
                                  6 & &
                \end{tabular}
                &
                \begin{minipage}[c]{68pt}
                    \begin{flushleft}
                    Form a new pile with \textbf{5}:
                    \end{flushleft}
                \end{minipage}
                &
                \begin{tabular}{l l l}
                       & & \\
                    4 & & \\
                    6 & \textbf{5} &
                \end{tabular}
            \end{tabular}\bigskip\bigskip

            \begin{tabular}{l p{56pt} l p{56pt} l p{56pt}}
                \begin{minipage}[c]{68pt}
                    \begin{flushleft}
                    Add \textbf{1} to the left-most pile:
                    \end{flushleft}
                \end{minipage}
                &
                \begin{tabular}{l l l}
                    \textbf{1} & & \\
                               4 & & \\
                               6 & 5 &
                \end{tabular}
                &
                \begin{minipage}[c]{68pt}
                    \begin{flushleft}
                    Form a new pile with \textbf{8}:
                    \end{flushleft}
                \end{minipage}
                &
                \begin{tabular}{l l l}
                    1 & & \\
                    4 & & \\
                    6 & 5 & \textbf{8}
                \end{tabular}
                &
                \begin{minipage}[c]{68pt}
                    \begin{flushleft}
                    Then place \textbf{7} atop 8:
                    \end{flushleft}
                \end{minipage}
                &
                \begin{tabular}{l l l}
                    1 & & \\
                    4 & & \textbf{7} \\
                    6 & 5 & 8
                \end{tabular}
            \end{tabular}\bigskip\bigskip

            \begin{tabular}{l p{56pt} l p{56pt} l p{56pt}}
                \begin{minipage}[c]{68pt}
                    \begin{flushleft}
                    Add \textbf{2} to the middle pile:
                    \end{flushleft}
                \end{minipage}
                &
                \begin{tabular}{l l l}
                    1 & & \\
                    4 & \textbf{2} & 7 \\
                    6 & 5 & 8
                \end{tabular}
                &
                \begin{minipage}[c]{68pt}
                    \begin{flushleft}
                    Finally, place \textbf{3} atop 7:
                    \end{flushleft}
                \end{minipage}
                &
                \begin{tabular}{l l l}
                    1 & & \textbf{3} \\
                    4 & 2 & 7 \\
                    6 & 5 & 8
                \end{tabular}
                &
                \begin{minipage}[c]{68pt}
                    $\phantom{to force spacing}$
                \end{minipage}
                &
                $\phantom{to force spacing}$
            \end{tabular}\bigskip

        \end{center}

        \noindent Then, by reading up the columns of $R(\sigma)$ from
        left to right, $RPW(R(\sigma)) = 64152873 \in
        S_{8}(3\textrm{-}\bar{1}\textrm{-}42)$.

    \end{exgr}

    Given $\sigma \in \mathfrak{S}_{n}$, the formation of $R(\sigma)$
    can be viewed as an iterated, non-recursive form of the Schensted
    Insertion Algorithm for interposing values into the rows of a
    standard Young tableau (see \cite{refAD1999}).  In
    \cite{refBLFPSAC05} the authors augment the formation of
    $R(\sigma)$ so that the resulting extension of Patience Sorting
    becomes a full non-recursive analogue of the celebrated
    Robinson-Schensted-Knuth (or RSK) Correspondence.  As with RSK,
    this Extended Patience Sorting Algorithm (given as Algorithm
    \ref{alg:ExtendedPSalgorithm} in Section
    \ref{sec:Introduction:ExtendedPSAlgorithm} below) takes a simple
    idea --- that of placing cards into piles --- and uses it to build
    a bijection between elements of the symmetric group
    $\mathfrak{S}_{n}$ and certain pairs of combinatorial objects.  In
    the case of RSK, one uses the Schensted Insertion Algorithm to
    build a bijection with (unrestricted) pairs of standard Young
    tableau having the same shape (see \cite{refSagan2000}).  However,
    in the case of Patience Sorting, one achieves a bijection between
    permutations and (somewhat more restricted) pairs of pile
    configurations having the same shape.  We denote this latter
    bijection by $\sigma \stackrel{PS}{\longleftrightarrow}(R(\sigma),
    S(\sigma))$ and call $R(\sigma)$ (resp.  $S(\sigma)$) the
    \emph{insertion piles} (resp.  \emph{recording piles})
    corresponding to $\sigma$.  Collectively, we also call
    $(R(\sigma), S(\sigma))$ the \emph{stable pair} of pile
    configurations corresponding to $\sigma$ and characterize such
    pairs in~\cite{refBLFPSAC05} using a somewhat involved pattern
    avoidance condition on their reverse patience words.

    Barred (generalized) permutation patterns like
    $3\textrm{-}\bar{1}\textrm{-}42$ arise quite naturally when
    studying Patience Sorting.  We discuss and enumerate the avoidance
    classes for several related patterns in Section
    \ref{sec:BarredPatterns}.  Then, in Section
    \ref{sec:RestrictedPatienceSorting}, we examine properties of
    Patience Sorting under restricted input that can be characterized
    using such patterns.  One such characterization, discussed in
    Section
    \ref{sec:RestrictedPatienceSorting:RestrictedPermutations}, is for
    the crossings in the initial iteration of the Geometric Patience
    Sorting Algorithm given by the authors in \cite{refBLPP05}.  This
    geometric form for the Extended Patience Sorting Algorithm is
    naturally dual to Viennot's Geometric RSK (originally defined in
    \cite{refViennot1977}) and gives, among other things, a geometric
    interpretation for the stable pairs of
    $3\textrm{-}\bar{1}\textrm{-}42$-avoiding permutations
    corresponding to a permutation under Extended Patience Sorting.
    However, unlike Viennot's geometric form for RSK, the shadow lines
    in Geometric Patience Sorting are allowed to cross.  While a
    complete characterization for these crossings is given in
    \cite{refBLPP05} in terms of the pile configurations formed, this
    new result is the first step in providing a characterization for
    the permutations involved in terms of barred pattern avoidance.

    \bigskip

    We close this introduction by describing both the Extended and
    Geometric Patience Sorting Algorithms.  We also briefly review the
    notation of generalized permutation patterns.

    \subsection{Extended and Geometric Patience Sorting}
    \label{sec:Introduction:ExtendedPSAlgorithm}

        Mallows' original ``patience sorting procedure'' can be
        extended to a full bijection between the symmetric group
        $\mathfrak{S}_{n}$ and certain restricted pairs of pile
        configurations using the following algorithm (which was first
        introduced in \cite{refBLFPSAC05}):

        \begin{alg}[Extended Patience Sorting
        Algorithm]\label{alg:ExtendedPSalgorithm} Given $\sigma =
        c_{1} c_{2} \cdots c_{n} \in \mathfrak{S}_{n}$, inductively
        build \emph{insertion piles} $R(\sigma) = \{r_{1}, r_{2},
        \ldots, r_{m}\}$ and \emph{recording piles} $S(\sigma) =
        \{s_{1}, s_{2}, \ldots, s_{m}\}$ as follows:\medskip

            \begin{itemize}
                \item Place the first card $c_{1}$ from the deck into
                a pile $r_{1}$ by itself, and set $s_{1} = \{1\}$.\medskip

                \item For each remaining card $c_{i}$ ($i = 2, \ldots,
                n$), consider the cards $d_{1}, d_{2}, \ldots, d_{k}$ atop
                the piles $r_{1}, r_{2}, \ldots, r_{k}$ that have already
                been formed.\medskip

                    \begin{itemize}
                        \item If $c_{i} > \max\{d_{1}, d_{2}, \ldots,
                        d_{k}\}$, then put $c_{i}$ into a new pile
                        $r_{k+1}$ by itself and set
                        $s_{k+1} = \{i\}$.\medskip

                        \item Otherwise, find the left-most card
                        $d_{j}$ that is larger than $c_{i}$ and put
                        the card $c_{i}$ atop\\ pile $r_{j}$ while
                        simultaneously putting $i$ at the bottom of
                        pile $s_{j}$.\smallskip

                    \end{itemize}
            \end{itemize}
        \end{alg}

        \noindent Note that the pile configurations that comprise a
        resulting stable pair must have the same ``shape'', which we
        define as follows:

        \begin{defn}
            Given a pile configuration $R=\{r_{1}, r_{2}, \ldots,
            r_{m}\}$ on $n$ cards, we call the composition $\gamma =
            (|r_{1}|, |r_{2}|, \ldots, |r_{m}|)$ of $n$ the
            \emph{shape} of $R$ and denote this by $\mathrm{sh}(R) =
            \gamma \models n$.
        \end{defn}

        The idea behind Algorithm~\ref{alg:ExtendedPSalgorithm} is
        that we are using the auxiliary pile configuration $S(\sigma)$
        to implicitly label the order in which the elements of the
        permutation $\sigma \in \mathfrak{S}_{n}$ are added to the
        usual Patience Sorting pile configuration $R(\sigma)$ (which
        we now call the ``insertion piles'' of $\sigma$ in this
        context by analogy to RSK).  It is clear that this information
        then allows us to uniquely reconstruct $\sigma$ by reversing
        the order in which the cards were played.  As with normal
        Patience Sorting, we visualize the pile configurations
        $R(\sigma)$ and $S(\sigma)$ by listing their constituent piles
        vertically as illustrated in the following example.

        \begin{exgr}\label{eg:ExtendedPSexample}
            Given $\sigma = 6 4 5 1 8 7 2 3 \in \mathfrak{S}_{8}$ from
            Example \ref{eg:NormalPSexample} above, we simultaneously
            form the following pile configurations with shape
            $\mathrm{sh}(R(\sigma)) = \mathrm{sh}(S(\sigma)) = (3, 2,
            3)$ under Extended Patience Sorting
            (Algorithm~\ref{alg:ExtendedPSalgorithm}):\\

            \begin{center}
                    \begin{tabular}{llcll}
                        \begin{minipage}[c]{24pt}
                            $R(\sigma) = $
                        \end{minipage}
                        &
                        \begin{tabular}{l l l}
                            1 &    & 3 \\
                            4 & 2 & 7 \\
                            6 & 5 & 8
                        \end{tabular}
                        &
                        \begin{minipage}[c]{18pt}
                            and
                        \end{minipage}
                        &
                        \begin{minipage}[c]{24pt}
                            $S(\sigma) = $
                        \end{minipage}
                        &
                        \begin{tabular}{l l l}
                           1 &   & 5 \\
                           2 & 3 & 6 \\
                           4 & 7 & 8
                        \end{tabular}
                    \end{tabular}
            \end{center}

            \bigskip

            \noindent Note that the insertion piles $R(\sigma)$ are
            the same as the pile configuration formed in Example
            \ref{eg:NormalPSexample} and that $RPW(S(6 4 5 1 8 7 2 3))
            = 42173865 \in S_{8}(3\textrm{-}\bar{1}\textrm{-}42)$.

        \end{exgr}

        \begin{figure}[t]
            \centering

            \begin{tabular}{ccc}

                \begin{minipage}[c]{2.125in}

                    \psset{xunit=0.125in,yunit=0.125in}

                    \begin{pspicture}(0,0)(9,9)

                        \psaxes{->}(9,9)

                        \pspolygon[fillcolor=lightgray, linecolor=gray,
                        fillstyle=solid](0,6)(1,6)(1,4)(2,4)(2,1)(4,1)(4,0)(0,0)

                        \psline[linecolor=darkgray,
                        linewidth=1pt](0,6)(1,6)(1,4)(2,4)(2,1)(4,1)(4,0)

                         \rput(1,6){{\Large $\bullet$}}%
                         \rput(2,4){{\Large $\bullet$}}%
                         \rput(3,5){{\Large $\bullet$}}%
                         \rput(4,1){{\Large $\bullet$}}%
                         \rput(5,8){{\Large $\bullet$}}%
                         \rput(6,7){{\Large $\bullet$}}%
                         \rput(7,2){{\Large $\bullet$}}%
                         \rput(8,3){{\Large $\bullet$}}%

                    \end{pspicture}\\[0.15in]

                      (a) Shadowline \\ $L^{(0)}_{1}(64518723)$.

                \end{minipage}
                &

                \begin{minipage}[c]{2.125in}

                    \psset{xunit=0.125in,yunit=0.125in}

                    \begin{pspicture}(0,0)(9,9)

                        \psaxes{->}(9,9)

                        \pspolygon[fillcolor=lightgray, linecolor=gray,
                        fillstyle=solid](0,5)(3,5)(3,2)(7,2)(7,0)(0,0)

                        \psline[linecolor=darkgray,
                        linewidth=1pt](0,5)(3,5)(3,2)(7,2)(7,0)

                         \rput(1,6){{\Large $\bullet$}}%
                         \rput(2,4){{\Large $\bullet$}}%
                         \rput(3,5){{\Large $\bullet$}}%
                         \rput(4,1){{\Large $\bullet$}}%
                         \rput(5,8){{\Large $\bullet$}}%
                         \rput(6,7){{\Large $\bullet$}}%
                         \rput(7,2){{\Large $\bullet$}}%
                         \rput(8,3){{\Large $\bullet$}}%

                    \end{pspicture}\\[0.15in]

                      (b) Shadowline\\ $L^{(0)}_{2}(64518723)$.

                \end{minipage}
                &

                \begin{minipage}[c]{2.125in}

                    \psset{xunit=0.125in,yunit=0.125in}

                    \begin{pspicture}(0,0)(9,9)

                        \psaxes{->}(9,9)

                        \pspolygon[fillcolor=lightgray, linecolor=gray,
                        fillstyle=solid](0,8)(5,8)(5,7)(6,7)(6,3)(8,3)(8,0)(0,0)

                        \psline[linecolor=darkgray,
                        linewidth=1pt](0,8)(5,8)(5,7)(6,7)(6,3)(8,3)(8,0)

                         \rput(1,6){{\Large $\bullet$}}%
                         \rput(2,4){{\Large $\bullet$}}%
                         \rput(3,5){{\Large $\bullet$}}%
                         \rput(4,1){{\Large $\bullet$}}%
                         \rput(5,8){{\Large $\bullet$}}%
                         \rput(6,7){{\Large $\bullet$}}%
                         \rput(7,2){{\Large $\bullet$}}%
                         \rput(8,3){{\Large $\bullet$}}%

                    \end{pspicture}\\[0.15in]

                      (c) Shadowline\\ $L^{(0)}_{3}(64518723)$.

                \end{minipage}\bigskip

                \\

                \begin{minipage}[c]{2.125in}

                    \psset{xunit=0.125in,yunit=0.125in}

                    \begin{pspicture}(0,0)(9,9)

                        \psline[linecolor=darkgray,
                        linewidth=1pt](0,5)(3,5)(3,2)(7,2)(7,0)

                        \psaxes{->}(9,9)

                        \psline[linecolor=darkgray,
                        linewidth=1pt](0,6)(1,6)(1,4)(2,4)(2,1)(4,1)(4,0)

                        \psline[linecolor=darkgray,
                        linewidth=1pt](0,5)(3,5)(3,2)(7,2)(7,0)

                        \psline[linecolor=darkgray,
                        linewidth=1pt](0,8)(5,8)(5,7)(6,7)(6,3)(8,3)(8,0)

                        \rput(1,6){{\Large $\bullet$}}%
                        \rput(2,4){{\Large $\bullet$}}%
                        \rput(3,5){{\Large $\bullet$}}%
                        \rput(4,1){{\Large $\bullet$}}%
                        \rput(5,8){{\Large $\bullet$}}%
                        \rput(6,7){{\Large $\bullet$}}%
                        \rput(7,2){{\Large $\bullet$}}%
                        \rput(8,3){{\Large $\bullet$}}%

                        \rput(1,4){{\Large $\odot$}}%
                        \rput(2,1){{\Large $\odot$}}%
                        \rput(3,2){{\Large $\odot$}}%
                        \rput(5,7){{\Large $\odot$}}%
                        \rput(6,3){{\Large $\odot$}}%

                    \end{pspicture}\\[0.15in]

                      (d) Salient points $(\odot)$ \\ for
                      $D^{(0)}(64518723)$.

                \end{minipage}
                &

                \begin{minipage}[c]{2.125in}

                    \psset{xunit=0.125in,yunit=0.125in}

                    \begin{pspicture}(0,0)(9,9)

                        \psaxes{->}(9,9)

                        \psline[linecolor=darkgray,
                        linewidth=1pt](0,4)(1,4)(1,1)(2,1)(2,0)

                        \psline[linecolor=darkgray,
                        linewidth=1pt](0,2)(3,2)(3,0)

                        \psline[linecolor=darkgray,
                        linewidth=1pt](0,7)(5,7)(5,3)(6,3)(6,0)

                         \rput(1,4){{\Large $\bullet$}}%
                         \rput(2,1){{\Large $\bullet$}}%
                         \rput(3,2){{\Large $\bullet$}}%
                         \rput(5,7){{\Large $\bullet$}}%
                         \rput(6,3){{\Large $\bullet$}}%

                    \end{pspicture}\\[0.15in]

                      (e) Shadow Diagram\\ $D^{(1)}(64518723)$.

                \end{minipage}
                &

                \begin{minipage}[c]{2.125in}

                    \psset{xunit=0.125in,yunit=0.125in}

                    \begin{pspicture}(0,0)(9,9)

                        \psaxes{->}(9,9)

                        \psline[linecolor=darkgray,
                        linewidth=1pt](0,1)(1,1)(1,0)

                        \psline[linecolor=darkgray,
                        linewidth=1pt](0,3)(5,3)(5,0)

                         \rput(1,1){{\Large $\bullet$}}%
                         \rput(5,3){{\Large $\bullet$}}%

                    \end{pspicture}\\[0.15in]

                      (f) Shadow Diagram\\ $D^{(2)}(64518723)$.

                \end{minipage}

            \end{tabular}\\

            \caption{Examples of Shadowline and Shadow Diagram
            Construction.}
            \label{fig:ShadowExample}
        \end{figure}
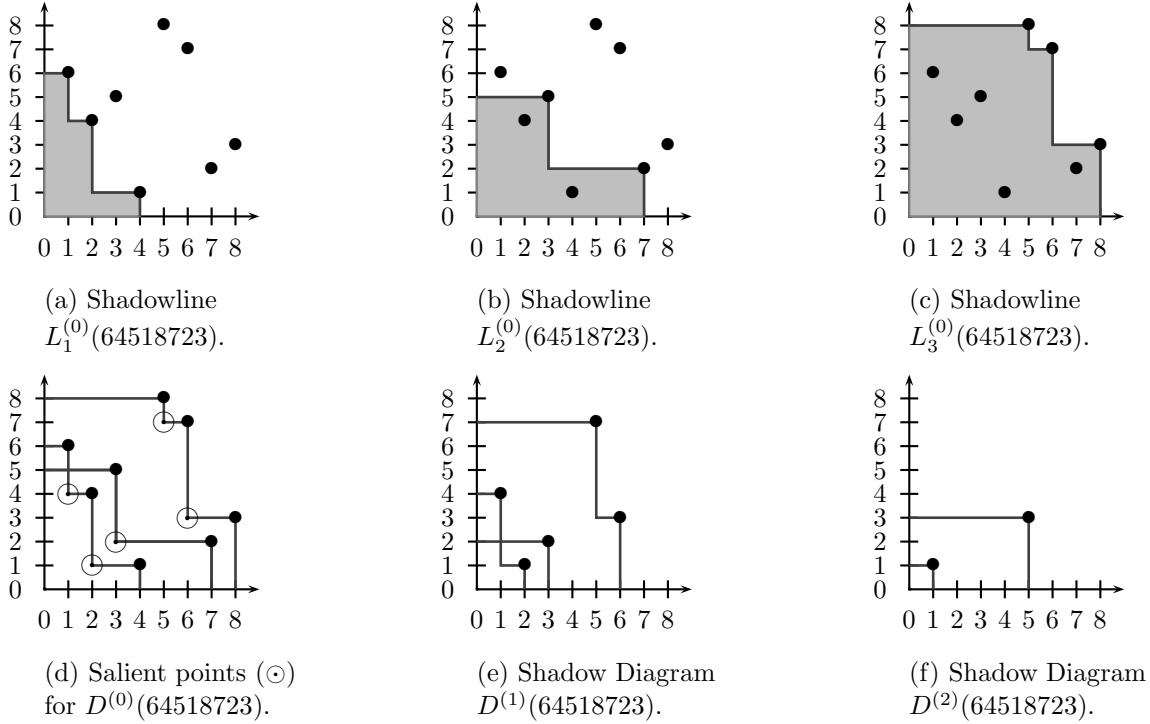

        In order to now describe a natural geometric form for this
        Extended Patience Sorting Algorithm, we begin with the
        following fundamental definition.

        \begin{defn}\label{defn:PSshadow}
            Given a lattice point $(m, n) \in \mathbb{Z}^{2}$, we
            define the \emph{(southwest) shadow} of $(m, n)$ to be the
            quarter space $U(m, n) = \{ (x, y) \in \mathbb{R}^{2}
            \ | \ x \leq m, \ y \leq n\}$.
        \end{defn}

        As with the northeasterly-oriented shadows that Viennot used
        when building his geometric form for RSK (see
        \cite{refViennot1977}), the most important use of these
        southwesterly-oriented shadows is in building shadowlines
        (which is illustrated in Figure
        \ref{fig:ShadowExample}(a)):

        \begin{defn}\label{defn:PSshadowline}
            The \emph{(southwest) shadowline} of $(m_{1}, n_{1}),
            (m_{2}, n_{2}), \ldots, (m_{k}, n_{k}) \in \mathbb{Z}^{2}$
            is defined to be the boundary of the union of the shadows
            $U(m_{1}, n_{1}), U(m_{2}, n_{2}), \ldots, U(m_{k},
            n_{k})$.
        \end{defn}

        In particular, we wish to associate to each permutation a
        certain collection of (southwest) shadowlines called its
        \emph{shadow diagram}.  However, unlike the
        northeasterly-oriented shadowlines used to define the
        northeast shadow diagrams of Geometric RSK
        \cite{refViennot1977}, these southwest shadowlines are allowed
        to intersect as illustrated in
        Figure~\ref{fig:ShadowExample}(d)--(e).  (We characterize
        those permutations having intersecting shadowlines under
        Definition \ref{defn:PSshadowDiagram} in Theorem
        \ref{thm:CharacterizedZerothCrossings} below.)

        \begin{defn}\label{defn:PSshadowDiagram}
            The \emph{(southwest) shadow diagram} $D^{(0)}(\sigma)$ of
            $\sigma = \sigma_{1}\sigma_{2}\cdots\sigma_{n} \in
            \mathfrak{S}_{n}$ consists of the (southwest) shadowlines
            $D^{(0)}(\sigma) = \{ L^{(0)}_{1}(\sigma),
            L^{(0)}_{2}(\sigma), \ldots, L^{(0)}_{k}(\sigma) \}$
            formed as follows:\medskip

            \begin{itemize}
                \item $L^{(0)}_{1}(\sigma)$ is the shadowline for
                those lattice points $(x, y) \in \{(1, \sigma_{1}),
                (2, \sigma_{2}), \ldots, (n, \sigma_{n})\}$ such that
                the shadow $U(x, y)$ does not contain any other
                lattice points.\medskip

                \item While at least one of the points $(1,
                \sigma_{1}), (2, \sigma_{2}), \ldots, (n, \sigma_{n})$
                is not contained in the shadowlines
                $L^{(0)}_{1}(\sigma), L^{(0)}_{2}(\sigma), \ldots,
                L^{(0)}_{j}(\sigma)$, define $L^{(0)}_{j+1}(\sigma)$
                to be the shadowline for the points
                \[
                (x, y) \in A := \{(i, \sigma_{i}) \ | \ (i,
                \sigma_{i}) \notin \bigcup^{j}_{k=1} L^{(0)}_{k}(\sigma)\}
                \]
                such that the shadow $U(x, y)$ does not contain any
                other lattice points from the set $A$.
            \end{itemize}

        \end{defn}

        \noindent In other words, we define a shadow diagram by
        inductively eliminating points in the permutation diagram
        until every point has been used to define a shadowline (as
        illustrated in Figure \ref{fig:ShadowExample}(a)--(c)).

        \medskip

        One can prove (see \cite{refBLFPSAC05}) that the ordinates
        (i.e., $y$-coordinates) of the points used to define each
        shadowline in the shadow diagram $D^{(0)}(\sigma)$ are exactly
        the \emph{left-to-right minima subsequences} (a.k.a.
        \emph{basic subsequences}) in the permutation $\sigma \in
        \mathfrak{S}_{n}$.  These are defined as follows:

        \begin{defn}\label{defn:LtoRminimaSubsequence}
            Let $\pi = \pi_{1} \pi_{2} \cdots \pi_{l}$ be a partial
            permutation on the set $[n] = \{1, 2, \ldots, n\}$.  Then
            the \emph{left-to-right minima} (resp.\ \emph{maxima})
            \emph{subsequence} of $\pi$ consists of those components
            $\pi_{j}$ of $\pi$ such that $\pi_{j} = \min\{ \pi_{i} \ |
            \ 1 \leq i \leq j\}$ (resp.  $\pi_{j} = \max\{ \pi_{i} \ |
            \ 1 \leq i \leq j\}$).
        \end{defn}

        \noindent We then inductively define the left-to-right minima
        (resp.\ maxima) subsequences $s_{1}, s_{2}, \ldots, s_{k}$ of
        the permutation $\sigma$ by taking $s_{1}$ to be the
        left-to-right minima (resp.\ maxima) subsequence for $\sigma$
        itself and then each remaining subsequence $s_{i}$ to be the
        left-to-right minima (resp.\ maxima) subsequence for the
        partial permutation obtained by removing the elements of
        $s_{1}, s_{2}, \ldots, s_{i-1}$ from~$\sigma$.

        \medskip

        Finally, one can produce a sequence $D(\sigma) =
        (D^{(0)}(\sigma), \ D^{(1)}(\sigma), \ D^{(2)}(\sigma), \
        \ldots$) of shadow diagrams for a given permutation $\sigma
        \in \mathfrak{S}_{n}$ by recursively applying Definition
        \ref{defn:PSshadowDiagram} to the southwest corners (called
        \emph{salient points}) of a given set of shadowlines (as
        illustrated in Figure \ref{fig:ShadowExample}(d)--(f)).  The
        only difference is that, with each iteration, newly formed
        shadowlines can only connect salient points along the same
        pre-existing shadowline.  One can then uniquely reconstruct
        the pile configurations $R(\sigma)$ and $S(\sigma)$ from these
        shadowlines by taking their intersections with the $x$- and
        $y$-axes in a certain canonical order (as detailed in
        \cite{refBLPP05}).

        \begin{defn}

            We call $D^{(k)}(\sigma)$ the $k^{\rm th}$ \emph{iterate}
            of the \emph{exhaustive shadow diagram} $D(\sigma)$ for
            the permutation $\sigma \in \mathfrak{S}_{n}$.

        \end{defn}

        \subsection{Generalized Pattern Avoidance}
        \label{sec:Introduction:GeneralizedPatternAvoidance}

        We first recall the following definition:

        \begin{defn}
            Let $\sigma = \sigma_{1}\sigma_{2}\cdots\sigma_{n} \in
            \mathfrak{S}_{n}$ and $\pi \in \mathfrak{S}_{m}$ with $m
            \leq n$.  Then we say that $\sigma$ \emph{contains} the
            \emph{(classical) permutation pattern} $\pi$ if there
            exists a subsequence $(\sigma_{i_{1}}, \sigma_{i_{2}},
            \ldots, \sigma_{i_{m}})$ of $\sigma$ (meaning
            $i_{1}~<~i_{2}~<~\cdots~<~i_{m}$) such that the word
            $\sigma_{i_{1}}\sigma_{i_{2}}\ldots\sigma_{i_{m}}$ is
            order-isomorphic to $\pi$.  I.e., each $\sigma_{i_{j}} <
            \sigma_{i_{j+1}}$ if and only if $\pi_{j} < \pi_{j+1}$.
        \end{defn}

        Note, though, that the elements in the subsequence
        $(\sigma_{i_{1}}, \sigma_{i_{2}}, \ldots, \sigma_{i_{m}})$ are
        not required to be contiguous in $\sigma$.  This motivates the

        \begin{defn}
            A \emph{generalized permutation pattern} is a classical
            permutation pattern $\pi$ in which one assumes that every
            element in the subsequence $(\sigma_{i_{1}},
            \sigma_{i_{2}}, \ldots, \sigma_{i_{m}})$ of $\sigma$ must
            be taken contiguously unless a dash is inserted between
            the corresponding order-isomorphic elements of the pattern
            $\pi$.
        \end{defn}

        Finally, if $\sigma$ does not contain a subsequence that is
        order-isomorphic to $\pi$, then we say that $\sigma$
        \emph{avoids} the pattern $\pi$. This motivated the

        \begin{defn}
            Given any collection $\pi^{(1)}, \pi^{(2)}, \ldots,
            \pi^{(k)}$ of permutation patterns (classical or
            generalized), we denote by $$S_{n}(\pi^{(1)}, \pi^{(2)},
            \ldots, \pi^{(k)}) = \bigcap_{i=1}^{k}S_{n}(\pi^{(i)}) =
            \bigcap_{i=1}^{k} \{\sigma \in \mathfrak{S}_{n} \ | \
            \sigma \textrm{ avoids } \pi^{(i)}\}$$
            the \emph{avoidance set} of permutations $\sigma \in
            \mathfrak{S}_{n}$ such that $\sigma$ simultaneously avoids
            each of the patterns $\pi^{(1)}, \pi^{(2)}, \ldots,\\
            \pi^{(k)}$.  Furthermore, the set $$\bigcup_{n \geq 1}
            S_{n}(\pi^{(1)}, \pi^{(2)}, \ldots, \pi^{(k)})$$ is called
            the \emph{(pattern) avoidance class} with \emph{basis}
            $\{\pi^{(1)}, \pi^{(2)}, \ldots, \pi^{(k)}\}$.
        \end{defn}

        More information about permutation patterns in general can be
        found in \cite{refBona2004}.

    \section{Barred and Unbarred Generalized Pattern Avoidance}
    \label{sec:BarredPatterns}

        An important further generalization of the notion of
        generalized permutation pattern requires that the context in
        which the occurrence of a generalized pattern occurs be taken
        into account.  The resulting concept of \emph{barred
        permutation patterns}, along with the accompanying notation,
        first arose within the study of stack-sortability of
        permutations by J. West \cite{refWest1990}.  Given how
        naturally these barred patterns now arise in the study of
        Patience Sorting (as illustrated in both \cite{refBLFPSAC05}
        and Section \ref{sec:RestrictedPatienceSorting} below), we
        initiate their systematic study in this section.

        \begin{defn}
            A \emph{barred (generalized) permutation pattern} $\beta$
            is a generalized permutation pattern in which overbars are
            used to indicate that barred values cannot occur at the
            barred positions.  As before, we denote by
            $S_{n}(\beta^{(1)}, \ldots, \beta^{(k)})$ the set of all
            permutations $\sigma \in \mathfrak{S}_{n}$ that
            simultaneously avoid $\beta^{(1)}, \ldots, \beta^{(k)}$
            (i.e., permutations that contain no subsequence that is
            order-isomorphic to any of the $\beta^{(1)}, \ldots,
            \beta^{(k)}$).
        \end{defn}

        \begin{exgr}
            A permutation $\sigma =
            \sigma_{1}\sigma_{2}\cdots\sigma_{n} \in \mathfrak{S}_{n}
            \notin
            S_{n}(3\textrm{-}\bar{5}\textrm{-}2\textrm{-}4\textrm{-}1)$
            contains an occurrence of the barred permutation pattern
            $3\textrm{-}\bar{5}\textrm{-}2\textrm{-}4\textrm{-}1$ if
            it contains an occurrence of the generalized pattern
            $3\textrm{-}2\textrm{-}4\textrm{-}1$ (i.e., contains a
            subsequence $(\sigma_{i_{1}}, \sigma_{i_{2}}, \ldots,
            \sigma_{i_{m}})$ that is order-isomorphic to the classical
            pattern 3241) in which no value larger than the element
            playing the role of ``4'' is allowed to occur between the
            elements playing the roles of ``3'' and ``2''.  This is
            one of the two basis elements for the pattern avoidance
            class used to characterize the set of 2-stack-sortable
            permutations \cite{refDGG, refDGW, refWest1990}.  (The
            other pattern is $2\textrm{-}3\textrm{-}4\textrm{-}1$,
            i.e., the classical pattern 2341.)  \end{exgr}\smallskip

        Despite the added complexity involved in avoiding barred
        permutation patterns, it is still sometimes possible to
        characterize the avoidance class for a barred permutation
        pattern in terms of an unbarred generalized permutation
        pattern.  The following theorem gives such a characterization
        for the pattern $3\textrm{-}\bar{1}\textrm{-}42$.  (Note,
        though, that there is no equivalent characterization for such
        barred permutation patterns as $1\bar{3}\textrm{-}42$ and
        $3\textrm{-}\bar{5}\textrm{-}2\textrm{-}4\textrm{-}1$.)

        \begin{thm}
        \label{thm:23-1}
        Let $B_{n} = \frac{1}{e}\sum_{k\geq 0}\frac{k^{n}}{k!}$ denote
        the $n^{\rm th}$ Bell number.  Then
            \begin{enumerate}
                \item $S_n(3\textrm{-}\bar{1}\textrm{-}42) =
                S_n(3\textrm{-}\bar{1}\textrm{-}4\textrm{-}2) =
                S_n(23\textrm{-}1)$

                \item $|S_n(3\textrm{-}\bar{1}\textrm{-}42)| = B_n$
            \end{enumerate}
        \end{thm}

        \begin{proof} (Sketch)\\
            \indent As in \cite{refClaesson2001}, we see that each of
            these sets consists of permutations having the form
            \[
            \sigma=\sigma_1 a_1\sigma_2a_2\dots\sigma_k a_k,
            \]
            where $a_k>a_{k-1}>\dots>a_2>a_1$ are the successive
            right-to-left minima of $\sigma$ (reversing the order of
            the elements in Definition
            \ref{defn:LtoRminimaSubsequence}) and where each segment
            $\sigma_i a_i$ is a decreasing subsequence.
        \end{proof}

        \medskip

        \begin{rem} \label{rem:why-3-1-42}

            We emphasize the following important consequences of
            Theorem \ref{thm:23-1}.\smallskip

            \begin{enumerate}

                \item Even though $S_n(3\textrm{-}\bar{1}\textrm{-}42)
                = S_n(23\textrm{-}1)$ by Theorem \ref{thm:23-1}(1), it
                is more natural to use avoidance of the barred pattern
                $3\textrm{-}\bar{1}\textrm{-}42$ in studying Patience
                Sorting.  As shown in \cite{refBLFPSAC05} and
                elaborated upon in Section
                \ref{sec:RestrictedPatienceSorting} below,
                $S_n(3\textrm{-}\bar{1}\textrm{-}42)$ is the set of
                equivalence classes of $\mathfrak{S}_n$ modulo the
                transitive closure of the relation
                $3\textrm{-}\bar{1}\textrm{-}42 \sim
                3\textrm{-}\bar{1}\textrm{-}24$.  (I.e., two
                permutations $\sigma, \tau \in \mathfrak{S}_{n}$ are
                equivalent if the elements creating an occurrence of
                one of these patterns in $\sigma$ form an occurrence
                of the other pattern in $\tau$.)  Moreover, each
                permutation $\sigma \in \mathfrak{S}_{n}$ in a given
                equivalence class has the same pile configuration
                $R(\sigma)$ under Patience Sorting, a description of
                which is significantly more difficult to describe for
                occurrences of the unbarred generalized permutation
                pattern $23\textrm{-}1$.

                \bigskip

                \item Marcus and Tardos proved in \cite{refMT2004}
                that the avoidance set $S_{n}(\pi)$ for any classical
                pattern $\pi$ grows at most exponentially fast as $n
                \to \infty$.  (This was previously known as the
                Stanley-Wilf Conjecture.)  The Bell numbers, though,
                satisfy $\log B_{n} = n (\log n - \log\log n + O(1))$
                and so exhibit superexponential growth.  (See
                \cite{refStanley1999} for more information about Bell
                numbers.)  While it was previously known that the
                Stanley-Wilf Conjecture does not extend to generalized
                permutation patterns (see, e.g.,
                \cite{refClaesson2001}), it took Theorem
                \ref{thm:23-1}(2) (originally proven in
                \cite{refBLFPSAC05} using Patience Sorting) to
                provide the first verification that one also cannot
                extend the Stanley-Wilf Conjecture to barred
                generalized permutation patterns.

                \bigskip

                \noindent A further abstraction of barred permutation
                pattern avoidance (called \emph{Bruhat-restricted
                avoidance}) was recently given by A.~Woo and A.~Yong
                in \cite{refWY2005}.  The result in Theorem
                \ref{thm:23-1}(2) has led A.~Woo to conjecture to the
                second author that the Stanley-Wilf ex-Conjecture also
                does not extend to this new notion of pattern
                avoidance.

            \end{enumerate}

        \end{rem}

        We conclude this section with a simple corollary to Theorem
        \ref{thm:23-1} that gives similar equivalences and
        enumerations for some barred permutation patterns that also
        arise naturally in the study of Patience Sorting (see
        Proposition \ref{prop:3-1-42} and Theorem
        \ref{thm:CharacterizedZerothCrossings} in Section
        \ref{sec:RestrictedPatienceSorting} below).

        \begin{cor} \label{cor:3-12}
            Using the notation in Theorem \ref{thm:23-1},

            \begin{enumerate}

                \item $S_n(31\textrm{-}\bar{4}\textrm{-}2) =
                S_n(3\textrm{-}1\textrm{-}\bar{4}\textrm{-}2) =
                S_n(3\textrm{-}12)$

                \item $S_n(\bar{2}\textrm{-}41\textrm{-}3) =
                S_n(\bar{2}\textrm{-}4\textrm{-}1\textrm{-}3) =
                S_n(2\textrm{-}4\textrm{-}1\textrm{-}\bar{3}) =
                S_n(2\textrm{-}41\textrm{-}\bar{3})$

                \item $|S_n(\bar{2}\textrm{-}41\textrm{-}3)| =
                |S_n(31\textrm{-}\bar{4}\textrm{-}2)| =
                |S_n(3\textrm{-}\bar{1}\textrm{-}42)| = B_n$.

            \end{enumerate}
        \end{cor}

        \begin{proof} (Sketches)
            \begin{enumerate}
                \item Take reverse complements in
                $S_n(3\textrm{-}\bar{1}\textrm{-}42)$ and apply Theorem
                \ref{thm:23-1}.

                \item Similar to (1).  (Note that (2) is also proven
                in \cite{refALR2005}.)

                \item This follows from the fact that the patterns
                $3\textrm{-}1\textrm{-}\bar{4}\textrm{-}2$ and
                $\bar{2}\textrm{-}4\textrm{-}1\textrm{-}3$ are
                inverses of each other.\\[-28pt]
            \end{enumerate}
        \end{proof}

        \smallskip

    \section{Patience Sorting under Restricted Input}
    \label{sec:RestrictedPatienceSorting}

        \subsection{Patience Sorting on Restricted Permutations}
        \label{sec:RestrictedPatienceSorting:RestrictedPermutations}

        The similarities between the Extended Patience Sorting
        Algorithm (Algorithm \ref{alg:ExtendedPSalgorithm}) and RSK
        applied to permutations is perhaps most observable in the
        following simple proposition:

        \begin{prop} \label{prop:monotone-patterns}
            Let $\textrm{\emph{\i}}_k =
            1\textrm{-}2\textrm{-}\cdots\textrm{-}k = 1 2 \cdots k$
            and $\textrm{\emph{\j}}_k =
            k\textrm{-}\cdots\textrm{-}2\textrm{-}1 =
            k \cdots 2 1$ be the classical monotone
            permutation patterns.  Then there is

            \begin{enumerate}

                \item a bijection between
                $S_n(\textrm{\emph{\i}}_{k+1})$ and pairs of pile
                configurations having the same composition shape
                $\gamma = (\gamma_{1}, \gamma_{2}, \ldots, \gamma_{m})
                \models n$ but with at most $k$
                piles (i.e., $m \leq k$).

                \item a bijection between
                $S_n(\textrm{\emph{\j}}_{k+1})$ and pairs of pile
                configurations having the same composition shape
                $\gamma = (\gamma_{1}, \gamma_{2}, \ldots, \gamma_{m})
                \models n$ but with no pile having
                more than $k$ cards in it (i.e., $\gamma_{i} \leq k$
                for each $i = 1, 2, \ldots, m$).

            \end{enumerate}
        \end{prop}

        \begin{proof} (Sketches)
            \begin{enumerate}
                \item Given $\sigma \in \mathfrak{S}_{n}$, a bijection
                is formed in \cite{refAD1999} between the set of piles
                $R(\sigma) = \{r_{1}, r_{2}, \ldots, r_{k}\}$ formed
                under Patience Sorting and the components of a
                particular distinguished longest increasing
                subsequence in $\sigma$.  Since avoiding the monotone
                pattern $\textrm{\i}_{k+1}$ is equivalent to
                restricting the length of the longest increasing
                subsequence in a permutation, the result then follows.

                \item Follows from (1) by reversing each of the
                permutations in $S_n(\textrm{\i}_{k+1})$ in
                order to form $S_n(\textrm{\j}_{k+1})$.\\[-28pt]
            \end{enumerate}
        \end{proof}

        Proposition \ref{prop:monotone-patterns} states that Patience
        Sorting can be used to efficiently compute the length of both
        the longest increasing and longest decreasing subsequences in
        a given permutation.  In particular, one can compute these
        lengths without needing to examine every subsequence of a
        permutation, just as with RSK. However, while both RSK and
        Patience Sorting can be used to implement this computation in
        $O(n\log(n))$ time, an extension of this technique is given in
        \cite{refBS2000} that also simultaneously tabulates all of the
        longest increasing or decreasing subsequences without
        incurring any additional asymptotic computational cost.

        \medskip

        As mentioned in Section \ref{sec:BarredPatterns} above,
        Patience Sorting also has immediate connections to certain
        barred permutation patterns:

        \begin{prop} \label{prop:3-1-42}

            $\phantom{to force a new line}$

            \begin{enumerate}
                \item $S_n(3\textrm{-}\bar{1}\textrm{-}42) = \{
                RPW(R(\sigma)) \mid \sigma \in \mathfrak{S}_{n} \}$.
                In particular, given $\sigma \in
                S_n(3\textrm{-}\bar{1}\textrm{-}42)$, the entries in
                each column of the insertion piles $R(\sigma)$ (when
                read from bottom to top) occupy successive positions
                in the permutation $\sigma$.

                \smallskip

                \item $S_n(\bar{2}\textrm{-}41\textrm{-}3) = \{
                RPW(R(\sigma^{-1})) \mid \sigma \in \mathfrak{S}_{n}
                \}$.  In particular, given $\sigma \in
                S_n(\bar{2}\textrm{-}41\textrm{-}3)$, the columns of
                the insertion piles $R(\sigma)$ (when read from top to
                bottom) contain successive values.
            \end{enumerate}

        \end{prop}

        \begin{proof}
            Part (1) is proven in \cite{refBLFPSAC05},
            and part (2) follows immediate by taking inverses in (1).
        \end{proof}

        As an immediate corollary, we can characterize an important
        category of classical permutation patterns in terms of barred
        permutation patterns.

        \begin{defn}
            Given a composition $\gamma = (\gamma_{1}, \gamma_{2},
            \ldots, \gamma_{m}) \models n$, the
            \emph{(classical) layered permutation pattern}
            $\pi_{\gamma} \in \mathfrak{S}_{n}$ is the permutation
            \[
                \gamma_{1} \cdot (\gamma_{1} - 1) \cdots 1 \cdot
                (\gamma_{1} + \gamma_{2}) \cdot (\gamma_{1} +
                \gamma_{2} - 1) \cdots (\gamma_{1} + 1) \cdots n \cdot
                (n - 1) \cdots (\gamma_{1} + \gamma_{2} + \cdots +
                \gamma_{m - 1} + 1).
            \]
        \end{defn}

        \medskip

        \begin{exgr}
            Given $\gamma = (3, 2, 3) \models 8$,
            the corresponding layered pattern is $\pi_{(3,2,3)} =
            \widehat{321} \widehat{54} \widehat{876} \in
            \mathfrak{S}_{8}$ (following the notation in
            \cite{refPrice1997}).  Moreover, applying Extended
            Patience Sorting (Algorithm~\ref{alg:ExtendedPSalgorithm})
            to $\pi_{(3,2,3)}$:\\

            \begin{center}
                    \begin{tabular}{llcll}
                        \begin{minipage}[c]{54pt}
                            $R(\pi_{(3,2,3)}) \ = $
                        \end{minipage}
                        &
                        \begin{tabular}{l l l}
                            1 &   & 6 \\
                            2 & 4 & 7 \\
                            3 & 5 & 8
                        \end{tabular}
                        &
                        \begin{minipage}[c]{18pt}
                            and
                        \end{minipage}
                        &
                        \begin{minipage}[c]{54pt}
                            $S(\pi_{(3,2,3)}) \ = $
                        \end{minipage}
                        &
                        \begin{tabular}{l l l}
                           1 &   & 6 \\
                           2 & 4 & 7 \\
                           3 & 5 & 8
                        \end{tabular}
                    \end{tabular}
            \end{center}

            \bigskip

            \noindent Note in particular that $\pi_{(3,2,3)}$
            satisfies both of the conditions in Proposition
            \ref{prop:3-1-42}, which illustrates the following
            characterization of layered patterns:

        \end{exgr}

        \begin{cor} \label{cor:layered}
           $S_n(3\textrm{-}\bar{1}\textrm{-}42,\bar{2}\textrm{-}41\textrm{-}3)$
           is the set of layered patterns in $\mathfrak{S}_{n}$.
        \end{cor}

        \begin{proof}
            Apply Proposition \ref{prop:3-1-42} noting that
            $S_n(3\textrm{-}\bar{1}\textrm{-}42,\bar{2}\textrm{-}41\textrm{-}3)
            = S_n(23\textrm{-}1,31\textrm{-}2)$ (as considered in
            \cite{refCM2002}).
        \end{proof}

        As a consequence of this interaction between Patience Sorting
        and barred permutation patterns, we can now explicitly
        characterize those permutations for which the initial
        iteration of Geometric Patience Sorting (as defined in Section
        \ref{sec:Introduction:ExtendedPSAlgorithm} above) yields
        non-crossing lattice paths.

        \begin{thm}
        \label{thm:CharacterizedZerothCrossings}
            The set $S_n(3\textrm{-}\bar{1}\textrm{-}42,
            31\textrm{-}\bar{4}\textrm{-}2)$ consists of all reverse
            patience words having non-intersecting shadow diagrams.
            (I.e., no shadowlines cross in the $0^{\rm th}$ iterate
            shadow diagram.)  Moreover, given a permutation $\sigma
            \in S_n(3\textrm{-}\bar{1}\textrm{-}42,
            31\textrm{-}\bar{4}\textrm{-}2)$, the values in the bottom
            rows of $R(\sigma)$ and $S(\sigma)$ increase from left to
            right.
        \end{thm}

        \begin{proof}
            From Theorem \ref{thm:23-1} and Corollary \ref{cor:3-12},
            $R(S_n(3\textrm{-}\bar{1}\textrm{-}42,
            31\textrm{-}\bar{4}\textrm{-}2)) =
            R(S_n(23\textrm{-}1,3\textrm{-}12))$ consists exactly of
            set partitions of $[n] = \{1, 2, \ldots, n\}$ whose
            components can be ordered so that both the minimal and
            maximal elements of the components simultaneously
            increase.  (These are called \emph{strongly monotone
            partitions} in \cite{refCMsubmitted}).

            Let $\sigma\in
            S_n(3\textrm{-}\bar{1}\textrm{-}42,31\textrm{-}\bar{4}\textrm{-}2)$.
            Since $\sigma$ avoids $3\textrm{-}\bar{1}\textrm{-}42$, we
            have that $\sigma=RPW(R(\sigma))$ by Proposition
            \ref{prop:3-1-42}.  Thus, the $i^{\rm th}$ shadowline
            $L^{(0)}_{i}(\sigma)$ of $\sigma$ is the boundary of the
            union of shadows with generating points in decreasing
            segments $\sigma_ia_i$, $i\in[k]$, where $\sigma_ia_i$ are
            as in the proof of Theorem \ref{thm:23-1}.  Let $b_i$ be
            the $i^{\rm th}$ left-to-right maximum of $\sigma$.  Then
            $b_i$ is the left-most (i.e. maximal) entry of
            $\sigma_ia_i$, so $\sigma_ia_i=b_i\sigma'_ia_i$ for some
            decreasing subsequence $\sigma'_i$.  Note that $\sigma'_i$
            may be empty so that $b_i=a_i$.

            Since $b_{i}$ is the $i^{\rm th}$ left-to-right maximum of
            $\sigma$, it must be at the bottom of the $i^{\rm th}$
            column of $R(\sigma)$ (similarly, $a_i$ is at the top of
            the $i^{\rm th}$ column).  So the bottom rows of both
            $R(\sigma)$ and $S(\sigma)$ must be in increasing order.

            Now consider the $i^{\rm th}$ and $j^{\rm th}$ shadowlines
            $L^{(0)}_{i}(\sigma)$ and $L^{(0)}_{j}(\sigma)$ of
            $\sigma$, respectively, where $i<j$.  We have that
            $b_i<b_j$ from which the initial horizontal segment of the
            $i^{\rm th}$ shadowline is lower than that of the $j^{\rm
            th}$ shadowline.  Moreover, $a_i$ is to the left of $b_j$,
            so the remaining segment of the $i^{\rm th}$ shadowline is
            completely to the left of the remaining segment of the
            $j^{\rm th}$ shadowline.  Thus, $L^{(0)}_{i}(\sigma)$ and
            $L^{(0)}_{j}(\sigma)$ do not intersect.
        \end{proof}

        In \cite{refBLPP05} the authors actually give the following
        stronger result:

        \begin{thm}\label{thm:NoncrossingPilesCondition}
            Each iterate $D_{SW}^{(m)}(\sigma)$ ($m\ge 0$) of $\sigma \in
            \mathfrak{S}_{n}$ is free from crossings if and only if
            every row in both $R(\sigma)$ and $S(\sigma)$ is monotone
            increasing from left to right.
        \end{thm}

        \noindent However, this only characterizes the output of the
        Extended Patience Sorting Algorithm involved.  As such,
        Theorem \ref{thm:CharacterizedZerothCrossings} provides the
        first step toward characterizing those permutations that
        result in non-crossing lattice paths under Geometric Patience
        Sorting.

        We conclude this section by noting that, while the strongly
        monotone condition implied by simultaneously avoiding
        $3\textrm{-}\bar{1}\textrm{-}42$ and
        $31\textrm{-}\bar{4}\textrm{-}2$ is necessary to alleviate
        such crossings, it is clearly not sufficient.  (The problem
        lies with what we call ``polygonal crossings'' in the shadow
        diagrams in \cite{refBLPP05}, which occur in permutations like
        $\sigma = 45312$.)  Thus, to avoid crossings at all iterations
        of Geometric Patience Sorting, we need to impose further
        ``ordinally increasing'' conditions on the set partition
        associated to a given permutation under Patience Sorting.  In
        particular, in addition to requiring just the minima and
        maxima elements in the set partition to increase as in the
        strongly monotone partitions encountered in the proof of
        Theorem \ref{thm:CharacterizedZerothCrossings}, it is
        necessary to require that every record value simultaneously
        increase under an appropriate ordering of the blocks.  That
        is, under a single ordering of these blocks, we must
        simultaneously have that the largest elements in each block
        increase, then the next largest elements, then the next-next
        largest elements, and so on.  E.g., the partition
        $\{\{5,3,1\},\{6,4,2\}\}$ of the set $[6] = \{1, 2, \ldots,
        6\}$ satisfies this condition.

        \subsection{Invertibility of Patience Sorting}
        \label{sec:RestrictedPatienceSorting:InvertibilityOfPS}

        It is clear that the pile configurations corresponding to two
        permutations under the Patience Sorting Algorithm need not be
        distinct in general (e.g., $R(3142) = R(3412)$).  As proven in
        \cite{refBLFPSAC05}, two permutations give rise to the same
        pile configuration under Patience Sorting if and only if they
        have the same left-to-right minima subsequences (e.g., $3142$
        and $3412$ both have the left-to-right minima subsequences
        $31$ and $42$).  In this section we characterize permutations
        having distinct pile configurations under Patience Sorting in
        terms of certain barred permutation patterns.  We then
        establish a non-trivial enumeration for the resulting
        avoidance sets.

        \begin{thm}\label{thm:unique}
            A pile configuration pile $R$ has a unique preimage
            $\sigma \in \mathfrak{S}_{n}$ under Patience Sorting if
            and only if $\sigma \in
            S_n(3\text{-}\bar{1}\text{-}42,3\text{-}\bar{1}\text{-}24)$.
        \end{thm}

        \begin{proof} (Sketch)

            It is clear that every pile configuration $R$ has at least
            one preimage, namely its reverse patience word $\sigma =
            RPW(R)$.  By Proposition \ref{prop:3-1-42}, reverse
            patience words are exactly those permutations that avoid
            the barred pattern $3\text{-}\bar{1}\text{-}42$.
            Furthermore, as shown in \cite{refBLFPSAC05}, two
            permutations have the same insertion piles under Extended
            Patience Sorting (Algorithm \ref{alg:ExtendedPSalgorithm})
            if and only if one can be obtained from the other by a
            sequence of order-isomorphic exchanges
            $3\text{-}\bar{1}\text{-}24\leadsto
            3\text{-}\bar{1}\text{-}42$ or
            $3\text{-}\bar{1}\text{-}42\leadsto
            3\text{-}\bar{1}\text{-}24$.  (I.e., the occurrence of one
            pattern is reordered to form an occurrence of the other
            pattern.)  Thus, it is easy to see that $R$ has a unique
            preimage $\sigma$ if and only if $\sigma$ has no
            occurrence of $3\text{-}\bar{1}\text{-}42$ or
            $3\text{-}\bar{1}\text{-}24$.
        \end{proof}

        \smallskip

        Given this pattern avoidance characterization of
        invertibility, we have the following recurrence relation for
        the number of permutations having distinct pile configurations
        under Patience Sorting:

        \begin{lem}
        \label{lem:RecurrenceRelation}
            Set $f(n) = |S_n(3\text{-}\bar{1}\text{-}42,
            3\text{-}\bar{1}\text{-}24)|$ and, for $k \le n$,
            \[
            f(n,k)=|\{\sigma \in
            S_n(3\text{-}\bar{1}\text{-}42,3\text{-}\bar{1}\text{-}24)
            \,:\,\sigma(1)=k \}|.
            \]
            Then $f(n)=\sum_{k=1}^{n}{f(n,k)}$, and we have the following
            recurrence relation for $f(n,k)$:
            \begin{eqnarray}
                f(n,0)=0 & \textrm{for} & n \geq 1
                \label{eq:RecurrenceRelation:Eqn1}\\
                f(n,1)=f(n,n)=f(n-1) & \textrm{for} & n \geq 1
                \label{eq:RecurrenceRelation:Eqn2}\\
                f(n,2)=0 & \textrm{for} & n \geq 3
                \label{eq:RecurrenceRelation:Eqn3}\\
                f(n,k)=f(n,k-1)+f(n-1,k-1)+f(n-2,k-2) & \textrm{for} &
                n \geq 3 \label{eq:RecurrenceRelation:Eqn4}
            \end{eqnarray}
            subject to the initial conditions $f(0) = f(0,0) = 1$.

        \end{lem}

        \begin{proof}
            First note that Equation
            \eqref{eq:RecurrenceRelation:Eqn1} is the obvious boundary
            condition for $k = 0$.

            Now suppose that the first letter of $\sigma \in
            S_n(3\text{-}\bar{1}\text{-}42,3\text{-}\bar{1}\text{-}24)$
            is $\sigma(1) = 1$ or $n$.  Then $\sigma(1)$ cannot be
            part of any occurrence of $3\text{-}\bar{1}\text{-}42$ or
            $3\text{-}\bar{1}\text{-}24$ in $\sigma$.  Thus, deletion
            of $\sigma(1)$, and subtraction of $1$ from each component
            if $\sigma(1) = 1$, yields a bijection with
            $S_{n-1}(3\text{-}\bar{1}\text{-}42,
            3\text{-}\bar{1}\text{-}24)$ so that Equation
            \eqref{eq:RecurrenceRelation:Eqn2} follows.

            Similarly, suppose that the first letter of $\sigma \in
            S_n(3\text{-}\bar{1}\text{-}42,3\text{-}\bar{1}\text{-}24)$
            is $\sigma(1) = 2$.  Then the first column of $R(\sigma)$
            must be {\tiny \begin{tabular}{c}1\\2 \end{tabular}}
            regardless of where 1 occurs in $\sigma$.  Therefore,
            $R(\sigma)$ has a unique preimage $\sigma$ if and only if
            $\sigma = 21 \in \mathfrak{S}_{2}$ so that Equation
            \eqref{eq:RecurrenceRelation:Eqn3} follows.

            Finally, suppose that $\sigma \in
            S_n(3\text{-}\bar{1}\text{-}42,3\text{-}\bar{1}\text{-}24)$
            with $3\le k\le n$.  Since $\sigma$ avoids
            $3\text{-}\bar{1}\text{-}42$, $\sigma$ is a RPW by
            Proposition \ref{prop:3-1-42}, and hence the left prefix
            of $\sigma$ from $k$ to $1$ is a decreasing subsequence.
            Let $\sigma'$ be the permutation obtained by interchanging
            the values $k$ and $k-1$ in $\sigma$.  Then the only
            instances of the patterns $3\text{-}\bar{1}\text{-}42$ and
            $3\text{-}\bar{1}\text{-}24$ in $\sigma'$ must involve
            both $k$ and $k-1$.  Note that the number of $\sigma$ for
            which no instances of these patterns are created by
            interchanging $k$ and $k-1$ is exactly $f(n,k-1)$.

            There are then two cases in which an instance of the
            barred pattern $3\text{-}\bar{1}\text{-}42$ or
            $3\text{-}\bar{1}\text{-}24$ will be created in $\sigma'$
            by this interchange:

            \emph{Case 1.} If $k-1$ occurs between $\sigma(1) = k$ and
            $1$ in $\sigma$, then $\sigma(2) = k - 1$, so
            interchanging $k$ and $k-1$ creates an instance of the
            pattern $23\text{-}1$ via the subsequence $(k-1,k,1)$ in
            $\sigma'$.  Thus, by Theorem \ref{thm:23-1}, $\sigma'$
            contains $3\text{-}\bar{1}\text{-}42$ from which $\sigma'
            \in S_{n}(3\text{-}\bar{1}\text{-}42)$ if and only if
            $k-1$ occurs after $1$ in $\sigma$.  Note also that if
            $\sigma(2) = k - 1$, then deleting $k$ yields a bijection
            with permutations in $S_{n-1}(3\text{-}\bar{1}\text{-}42,
            3\text{-}\bar{1}\text{-}24)$ that start with $k-1$.  So
            the number of permutations counted in Case 1 is exactly
            $f(n-1,k-1)$.

            \emph{Case 2.} If $k-1$ occurs to the right of $1$ in
            $\sigma$, then $\sigma'$ both contains the subsequence
            $(k-1,1,k)$ and avoids the pattern
            $3\text{-}\bar{1}\text{-}42$, so it must also contain the
            pattern $3\text{-}\bar{1}\text{-}24$.  If an instance of
            $3\text{-}\bar{1}\text{-}24$ in $\sigma'$ involves both
            $k-1$ and $k$, then $k-1$ and $k$ must play the roles of
            ``3'' and ``4'', respectively.  If the value $\ell$
            preceding $k$ is not $1$, then the subsequence
            $(k-1,1,\ell,k)$ is an instance of $3\text{-}1\text{-}24$,
            so $(k-1,\ell,k)$ is not an instance of
            $3\text{-}\bar{1}\text{-}24$.  Therefore, for $\sigma'$ to
            contain $3\text{-}\bar{1}\text{-}24$, $k$ must follow $1$
            in $\sigma'$, and so $k-1$ follows $1$ in $\sigma$.  If
            the letter preceding $1$ is some $m < k$, then the
            subsequence $(m,1,k-1)$ is an instance of
            $3\text{-}\bar{1}\text{-}24$ in $\sigma$, which is
            impossible.  Therefore, $k$ must precede $1$ in $\sigma$,
            from which $\sigma$ must start with the initial segment
            $(k,1,k-1)$.  But then deleting the values $k$ and $1$ and
            then subtracting $1$ from each component yields a
            bijection with permutations in
            $S_{n-2}(3\text{-}\bar{1}\text{-}42,3\text{-}\bar{1}\text{-}24)$
            that start with $k-2$.  It follows that the number of
            permutations counted in Case 2 is then exactly
            $f(n-2,k-2)$, which yields Equation
            \eqref{eq:RecurrenceRelation:Eqn4}.
        \end{proof}

        \medskip

        If we denote by
        \[
            \Phi(x, y) = \sum_{n = 0}^{\infty} \sum_{k=0}^{n}
            f(n,k)x^{n}y^{k}
        \]

        \medskip

        \noindent the bivariate generating
        function for the sequence $\{f(n,k)\}_{n \geq k\geq 0}$, then
        Equation \eqref{eq:RecurrenceRelation:Eqn4} implies that

        \[
            (1-y-xy-x^{2}y^{2})\Phi(x,y) =
            1-y-xy+xy^{2}-xy^{2}\Phi(xy,1)+xy(1-y-xy)\Phi(x,1).
        \]

        \medskip

        \noindent Moreover, using the kernel method, one can show that

        \[
            x + 1 + \frac{\sqrt{1+2x+5x^2}-x-1}{2} \cdot F(x) -
            F\left(\frac{\sqrt{1+2x+5x^2}-x-1}{2x}\right) = 0
        \]

        \medskip

        \noindent where $F(x) = \sum_{n \ge 0}f(n)x^{n}=\Phi(x,1)$ is the
        generating function for the sequence $\{f(n)\}_{n \ge 0}$.

        \medskip

        We conclude with the following main enumerative result about
        invertibility of Patience Sorting.

        \begin{thm}
        \label{thm:EnumerateInvertibility}
            Denote by $F_{n}$ the $n^{\rm th}$ Fibonacci number (with
            $F_{0} = F_{1} = 1$) and by

            \begin{displaymath}
                a(n,k) \hspace{0.25cm} = \hspace{-0.75cm}
                \sum_{\stackrel{n_1,\dots,n_k\ge
                0}{n_1+\dots+n_k=n-k-2}}\hspace{-0.75cm}{F_{n_1}F_{n_2}\dots
                F_{n_k}}
            \end{displaymath}

            \medskip

            convolved Fibonacci numbers for $n\ge k+2$ (where $a(n,k)
            := 0$ otherwise).  Then, defining

            \[
                X =
                    \begin{bmatrix}
                        f(0)  \\
                        f(1)  \\
                        f(2)  \\
                        f(3)  \\
                        f(4)  \\
                        \vdots
                    \end{bmatrix},
                \quad
                F =
                    \begin{bmatrix}
                        1  \\
                        F_{0} \\
                        F_{1} \\
                        F_{2} \\
                        F_{3} \\
                        \vdots
                    \end{bmatrix},
                \ \ \mathrm{and} \quad
                \mathbf{A} = (a(n,k))_{n,k\ge 0} =
                    \begin{bmatrix}
                        0 &  &  &  &  &  \\
                        0 & 0 &  &  &  &  \\
                        a(2,0) & 0 & 0 &  &  &  \\
                        a(3,0) & a(3,1) & 0 & 0 &  &  \\
                        a(4,0) & a(4,1) & a(4,2) & 0 & \phantom{aa}0\phantom{aa} &  \\
                       \vdots & \vdots & \vdots & \vdots & \vdots & \ddots \\
                    \end{bmatrix},
            \]

            \bigskip

            \indent we have that $X = (\mathbf{I} -
            \mathbf{A})^{-1}F$, where $\mathbf{I}$ is the infinite
            identity matrix and $\mathbf{A}$ is lower triangular.

        \end{thm}

        \begin{proof} (Sketch)

            From Equations
            \eqref{eq:RecurrenceRelation:Eqn1}--\eqref{eq:RecurrenceRelation:Eqn4},
            we can conjecture an equivalent recurrence where
            (\ref{eq:RecurrenceRelation:Eqn3}) and
            (\ref{eq:RecurrenceRelation:Eqn4}) are replaced by the
            following equation (here $\delta_{nk}$ is the Kronecker
            delta function):

            \begin{equation}
            \label{eq:RecurrenceRelation:Eqn5}
                f(n,k)=\sum_{m=0}^{k-3}{c(k,m)f(n-k+m)} +
                \delta_{nk}F_{k-2}, \quad n\ge k\ge 2.
            \end{equation}

            \noindent For this relation to hold, the coefficients
            $c(k,m)$ must satisfy the following recurrence relation:
            \[
                c(k,m) = c(k-1,m-1) + c(k-1,m) + c(k-2,m), \quad k\ge
                2,
            \]
            or, equivalently,
            \[
                c(k-1,m-1) = c(k,m) - c(k-1,m) - c(k-2,m), \quad k\ge
                2,
            \]

            \noindent with $c(2,0)=1$ and $c(k,m)=0$ in the case that
            $k < 2$, $m < 0$ or $m > k - 2$.  This implies that the
            generating function for the sequence $\{c(k,m)\}_{k\ge 0}$
            (for each $m\ge 0$) is

            \[
            \sum_{n\ge 0}{c(k,m)x^k} = \frac{x^{m+2}}{(1-x-x^2)^{m+1}}.
            \]

            \noindent It follows that the coefficients $c(k,m) =
            a(k,m)$ in Equation \eqref{eq:RecurrenceRelation:Eqn5} are
            convolved Fibonacci numbers \cite{refOEIS} forming the
            so-called skew Fibonacci-Pascal triangle in the matrix
            $\mathbf{A}=(a(k,m))_{k,m\ge 0}$.  In particular, the
            sequence of nonzero entries in column $m\ge 0$ of
            $\mathbf{A}$ is the $m^{\rm th}$ convolution of the
            sequence $\{F_n\}_{n\ge 0}$.

            Combining the expansion of $f(n,n)$ from Equation
            \eqref{eq:RecurrenceRelation:Eqn5} with Equation
            \eqref{eq:RecurrenceRelation:Eqn2}, we obtain

            \[
            f(n) = \sum_{m=0}^{n-2}{a(n,m)f(m)} + F_{n-1},
            \]

            \noindent which is equivalent to the matrix equation $X =
            \mathbf{A}X + F$.  Since $\mathbf{I} - \mathbf{A}$ is
            clearly invertible, the result follows.
        \end{proof}

        Due to space restrictions, we omit a direct bijective proof
        of Theorem \ref{thm:EnumerateInvertibility} that will be
        included in the full article.

        \begin{rem}
            Note that $\mathbf{A}$ is a strictly lower triangular
            matrix with zero sub-diagonal.  From this it follows that
            multiplication of a matrix $\mathbf{B}$ by $\mathbf{A}$
            shifts the position of the highest nonzero diagonal in
            $\mathbf{B}$ down by two rows, so $(\mathbf{I} -
            \mathbf{A})^{-1}=\sum_{n\ge 0}{\mathbf{A}^n}$ as a Neumann
            series, and thus all nonzero entries of $(\mathbf{I} -
            \mathbf{A})^{-1}$ are positive integers.
        \end{rem}

        Finally, one can explicitly compute

        \[
            \mathbf{A} =
                \begin{bmatrix}
                    0 \\
                    0 &  0 \\
                    1 &  0 &  0 \\
                    1 &  1 &  0 &  0 \\
                    2 &  2 &  1 &  0 &  0 \\
                    3 &  5 &  3 &  1 &  0 &  0 \\
                    5 & 10 &  9 &  4 &  1 &  0 &  0 \\
                    8 & 20 & 22 & 14 &  5 &  1 &  0 &  0 \\
                    \vdots &  \vdots &  \vdots &  \vdots & \vdots & \vdots &
                    \vdots & \vdots & \ddots \\
                \end{bmatrix}
                \quad \Longrightarrow \quad
            (\mathbf{I} - \mathbf{A})^{-1} =
                \begin{bmatrix}
                     1 \\
                     0 &  1 \\
                     1 &  0 &  1 \\
                     1 &  1 &  0 &  1 \\
                     3 &  2 &  1 &  0 &  1 \\
                     7 &  6 &  3 &  1 &  0 &  1 \\
                    21 & 16 & 10 &  4 &  1 &  0 &  1 \\
                    66 & 50 & 30 & 15 &  5 &  1 &  0 & 1 \\
                    \vdots &  \vdots &  \vdots &  \vdots &  \vdots & \vdots & \vdots & \vdots & \ddots\\
                \end{bmatrix}
        \]

        \bigskip

        \noindent from which the first few values of the sequence
        $\{f(n)\}_{n \ge 0}$ are immediately calculable as
        \[
        1, 1, 2, 4, 9, 23, 66, 209, 718,
        2645, 10373, 43090, 188803, 869191, 4189511, 21077302,
        110389321 \ldots\ .
        \]

\end{document}